\documentclass[12pt]{article}

\usepackage{amsmath,amssymb}
\usepackage{cite}
\usepackage{url}
\usepackage{placeins}
\usepackage{latexsym}
\usepackage{graphicx}
\usepackage{subfigure}
\usepackage{amssymb}
\usepackage{mathrsfs}
\usepackage{cases}
\usepackage{amsfonts}
\usepackage{indentfirst}
\usepackage{nopageno}
\usepackage{algorithmic}
\usepackage{algorithm}
\newcommand{\op}{\mathcal{O}}
\newcommand{\ie}{\emph{i.e.}}

\newcommand{\R}{\mathbb{R}}

\newtheorem{Proposition}{\bf  Proposition}
\newtheorem{Definition}{\bf  Definition}
\newtheorem{Corollary}{\bf  Corollary}

\begin{document}

\title{A class of linear solvers built on the Biconjugate $A$-Orthonormalization 
Procedure for solving unsymmetric linear systems}
\author{
Bruno Carpentieri\thanks{Email: b.carpentieri@rug.nl, bcarpentieri@gmail.com (B. Carpentieri)}\\
{\small{\it Institute of Mathematics and Computing Science,}}\\
{\small{\it University of Groningen}}\\
{\small{\it Nijenborgh 9, PO Box 407, 9700 AK Groningen,
Netherlands}}\\
Yan-Fei Jing\thanks{E-mail: yanfeijing@uestc.edu.cn, 00jyfvictory@163.com (Y.-F. Jing)},\mbox{ }
Ting-Zhu Huang\thanks{E-mail: tzhuang@uestc.edu.cn,
tingzhuhuang@126.com (T.-Z. Huang)} \mbox{ }
Yong Duan\thanks{E-mail: duanyong@yahoo.cn (Y. Duan)}\\
{\small{\it School of Mathematical Sciences/Institute of Computational Science,}}\\
{\small{\it University of Electronic Science and Technology of China,}}\\
{\small{\it Chengdu, Sichuan, 611731, P. R. China}}}

\bibliographystyle{plain}

 \maketitle

\begin{abstract}

We present economical iterative algorithms built on the Biconjugate $A$-Orthonormalization 
Procedure for real unsymmetric and complex non-Hermitian systems. The principal characteristics
of the developed solvers is that they are fast convergent and cheap in memory. 
We report on a large combination of numerical experiments to demonstrate that the 
proposed family of methods is highly competitive and often superior to other popular
algorithms built upon the Arnoldi method and the biconjugate Lanczos procedures 
for unsymmetric linear sytems. 

\smallskip


\emph{AMSC}: 65F10

\end{abstract}

\section{Introduction}~\label{sec:intro}

In this study we investigate variants of the Lanczos method
for the iterative solution of real unsymmetric and/or complex non-Hermitian 
linear systems

\begin{equation}\label{eq:linsys}
Ax=b, 
\end{equation}

with the motivation of obtaining smoother and, hopefully, faster convergence 
behavior in comparison with the BiCG method 
as well as its two evolving variants - the CGS method and one of the most 
popular methods in use today - the Biconjugate Gradient Stabilized (BiCGSTAB) 
method. 

Iterative methods for solving unsymmetric systems
are commonly developed upon the Arnoldi or the Lanczos biconjugate algorithms. 
These procedures generate an orthonormal basis for the Krylov 
subspaces associated with $A$ and an initial vector $v$, and
require only matrix-vector products by $A$. 
The generation of the vector recurrence by Arnoldi produces Hessenberg matrices, 
while the unsymmetric Lanczos biconjugation produces tridiagonal matrices. 
The price to pay due long recurrences in Arnoldi is the 
increasing orthogonalization cost along the iterations. 
In this paper we develop economical iterative algorithms built on the 
Biconjugate $A$-Orthonormalization Procedure 
presented in Section~\ref{sec:uno}.
The method ideally builds up a pair of Biconjugate $A$-Orthonormal 
(or, briefly, $A$-biorthonormal) basis
for the dual Krylov subspaces $K_m(A; v_1)$ 
and $A^T K_m(A^T;w_1)$ in the real case (which is $A^H K_m(A^H;w_1)$ in the complex case). The projection matrix
onto the corresponding Krylov subspace is tridiagonal, so that
the generation of the vector recurrences is extremely cheap in memory.
We provide the theoretical background for the developed algorithms and we discuss 
computational aspects.
We show by numerical experiments that the Biconjugate $A$-Orthonormalization Procedure
may lead to highly competitive solvers that are often superior
to other popular methods,~e.g. CGS, BiCGSTAB, IDR($s$).
We apply these techniques to sparse and dense matrix problems, 
both in real and complex arithmetic, 
arising from realistic applications in different areas. 
This study integrates and extends the preliminary investigations 
reported in~\cite{hgzl:09}, limited to the case of complex 
non-Hermitian systems.
In this paper we consider a much larger combination of experiments with 
both real and complex matrices having size order twice as large. 
Finally, we complete our work with a case study with 
dense linear systems in electromagnetic scattering from large structures.

The paper is structured as follows: in Section~\ref{sec:uno} we present 
the Biconjugate $A$-Orthonormalization Procedure and its properties; 
in Section~\ref{sec:due} we describe a general framework to derive linear 
solver from the proposed procedure, and we present the algorithmic
development of two Krylov projection algorithms. 
Finally, in Section~\ref{sec:tre} 
we report on extensive numerical experiments for solving 
large sparse and/or dense linear systems,
both real and complex.

\section{A general two-sided unsymmetric Lanczos biconjugate A-orthonormalization procedure}~\label{sec:uno}

Throughout this paper we denote the standard inner product of two real vectors $u, v \in \R^n$ as
\[
\left\langle {u,v} \right\rangle  = u^T v = \sum\limits_{i = 1}^n {u_i v_i }.
\]

Given two vectors $v_1$ and $w_1$ with euclidean inner product 
$\left\langle {\omega _1 ,Av_1 } \right\rangle  = 1$, we define Lanczos-type vectors 
$v_j$, $w_j$ and scalars $\delta_j$, $\beta_j$, $j = 1,2, \ldots ,m$ by the following recursions

\begin{align}
& \delta _{j + 1} v_{j + 1} = Av_j  - \beta _j v_{j - 1}  - \alpha _j v_j = s_{j+1}, & \label{rec:1} \\ 
& \beta _{j + 1} w_{j + 1} = A^T w_j  - \delta _j w_{j - 1}  - \alpha _j w_j  = t_{j+1}. & \label{rec:2} 
 \end{align}

where the scalars are chosen as

\[
 \alpha _j  = \left\langle {\omega _j ,A\left( {Av_j } \right)} \right\rangle,~  
 \displaystyle \delta _{j + 1}  = \left| {\left\langle {\hat \omega _{j + 1} ,A\hat v_{j + 1} } \right\rangle } \right|^{\frac{1}{2}},~ 
 \displaystyle \beta _{j + 1}  = \frac{\left\langle {\hat \omega _{j + 1} ,A\hat v_{j + 1}} \right\rangle}{\delta _{j + 1}}.
\]

This choice of the scalars guarantees that the recursions generate sequences of 
biconjugate $A$-orthonormal vector (or briefly, $A$-biorthonormal vectors) $v_i$ 
and $w_i$, according to the following definition

\begin{Definition}
Right and left Lanczos-type vectors $v_j,~j = 1, 2, \ldots, m$ and $w_i,~i = 1, 2, \ldots, m$ form 
a biconjugate A-orthonormal system in exact arithmetic, if and only if

\[
    \left\langle {\omega _i ,Av_j } \right\rangle  = \delta _{i,j} ,~1 \le i,j \le m.
\]
\end{Definition}

\begin{flushright}
 $\blacksquare$
\end{flushright}

Eqns.~(\ref{rec:1}-\ref{rec:2}) can be interpreted as a two-side Gram-Schmidt 
orthonormalization procedure where at step $i$ we multiply vectors $v_i$ and $w_i$ 
by $A$ and $A^T$, respectively, and we orthonormalize them against the most recently 
generated Lanczos-type pairs $(v_i, w_i)$ and $(v_{i-1}, w_{i-1})$. The vectors $\alpha_i v_i$, 
$\alpha_i w_i$ are the complex biconjugate $A$-orthonormal projections of $A v_i$ 
and $A^T w_i$ onto the most recently computed vectors $v_i$ and $w_i$; analogously, 
the vectors $\beta_i v_{i-1}$, $\delta_i w_{i-1}$ are the complex biconjugate 
$A$-orthonormalization projections of $A v_i$ and $A^T w_i$ onto the next computed 
vectors $v_{i-1}$ and $w_{i-1}$. The two sets of scalars satisfy the following relation 
\[
  \beta_{i+1} \delta_{i+1} = s_{i+1}^T A t_{i+1} = \left\langle {s_{i+1} ,A t_{i+1} } \right\rangle.
\]
The scalars $\beta_i$ and $\delta_i$ can be chosen with some freedom, provided the 
biconjugate $A$-orthonormalization property holds.
We sketch the complete procedure in~Algorithm~\ref{alg:biorth}

\begin{algorithm}
 \caption{The Lanczos A-biorthonormalization procedure}\label{alg:biorth}
\begin{algorithmic}[1]
  \STATE{\it Choose $v_1, \omega_1$, such that $\left\langle {\omega _1 ,Av_1 } \right\rangle  = 1$}
  \STATE Set $\beta _1  = \delta _1  \equiv 0,\omega _0  = v_0  = {\bf 0} \in \R^n $
  \FOR{$j=1,2,\ldots$}
     \STATE $\alpha _j  = \left\langle {\omega _j ,A\left( {Av_j } \right)} \right\rangle$
  \STATE $\hat v_{j + 1}  = Av_j  - \alpha _j v_j  - \beta _j v_{j - 1} $
  \STATE $\hat \omega _{j + 1}  = A^T \omega _j  - \alpha _j \omega _j  - \delta _j \omega _{j - 1} $
  \STATE $\displaystyle \delta _{j + 1}  = \left| {\left\langle {\hat \omega _{j + 1} ,A\hat v_{j + 1} } \right\rangle } \right|^{\frac{1}{2}} $
  \STATE $\displaystyle \beta _{j + 1}  = \frac{\left\langle {\hat \omega _{j + 1} ,A\hat v_{j + 1}} \right\rangle}{\delta _{j + 1}}$
  \STATE $\displaystyle v_{j + 1}  = \frac{{\hat v_{j + 1} }}{{\delta _{j + 1} }}$
  \STATE $\displaystyle \omega _{j + 1}  = \frac{{\hat \omega _{j + 1} }}{{\beta _{j + 1} }}$
\ENDFOR
\end{algorithmic}
\end{algorithm}

Notice that there is a clear analogy with the standard unsymmetric biconjugate Lanczos 
recursions. The matrix $A$ is not modified and is accessed only via matrix-vector 
products by $A$ and $A^T$. Similarly to the standard Lanczos procedure, the two most 
recently computed pairs of Lanczos-type vectors $v_k$ and $w_k$ for $k=i, i-1$ are 
needed at each step. These two vectors may be overwritten with the most recent updates. 
Therefore the memory storage is very limited compared to the Arnoldi method.
The price to pay is some lack of robustness due to possible vanishing of the inner 
products. Observe that the above algorithm is possible to breakdown whenever 
$\delta_{j+1}$ vanishes while $\hat w_{j+1}$ and $A \hat v_{j+1}$ are not equal to 
$0 \in \R^n$ appearing in line~7. In the interest 
of counteractions against such breakdowns, refer oneself to remedies such as so-called 
look-ahead strategies~\cite{patl:85,parl:92,frgn:93,gutk:97} which can enhance stability 
while increasing cost modestly, or others for example~\cite{day:97}. But that is outside 
the scope of this paper and we shall not pursue that here; for more details, please refer 
to~\cite{saad:96} and the references therein. In our experiments we never observed a 
breakdown of the algorithm, as we will see in Section~\ref{sec:tre}. However, it is fair to
mention that this problem may occurr.
The following proposition states some useful properties of Algorithm~\ref{alg:biorth}. 

\smallskip

\begin{Proposition}~\label{theo:biorth}
If Algorithm~\ref{alg:biorth} proceeds $m$ steps, then the right and left Lanczos-type vectors $v_j, j = 1, 2, \ldots, m$ and $w_i, i = 1, 2, \ldots, m$ form a biconjugate A-orthonormal system in exact arithmetic, i.e.,

\[
    \left\langle {\omega _i ,Av_j } \right\rangle  = \delta _{i,j} ,1 \le i,j \le m.
\]

\smallskip

Furthermore, denote by $V_m  = \left[ {v_1 ,v_2 , \ldots ,v_m } \right]$
and $W_m  = \left[ {w_1 ,w_2 , \ldots ,w_m } \right]$ the $n \times m$
matrices and by $\underline {T_m } $ the extended tridiagonal matrix of the form

\begin{equation}
\underline {T_m }  = \left[ \begin{array}{l}
 T_m  \\
 \delta _{m + 1} e_m^T  \\
 \end{array} \right], 
\end{equation}

where

\[
T_m  = \left[ {\begin{array}{*{20}c}
   {\alpha _1 } & {\beta _2 } & {} & {} & {}  \\
   {\delta _2 } & {\alpha _2 } & {\beta _3 } & {} & {}  \\
   {} &  \ddots  &  \ddots  &  \ddots  & {}  \\
   {} & {} & {\delta _{m - 1} } & {\alpha _{m - 1} } & {\beta _m }  \\
   {} & {} & {} & {\delta _m } & {\alpha _m }  \\
\end{array}} \right],
\]

whose entries are the coefficients generated during the algorithm implementation, and in which $\alpha _1 , \ldots ,\alpha _m ,\beta _2 , \ldots ,\beta _m $ are complex while $\delta _2 , \ldots ,\delta _m $ positive. 
Then with the Biconjugate A-Orthonormalization Procedure, the following four relations hold

\begin{align}
 & AV_m  = V_m T_m  + \delta _{m + 1} v_{m + 1} e_m^T & \label{eq:uno1}  \\
 & A^T W_m  = W_m T_m^T  + \beta _{m + 1} \omega _{m + 1} e_m^T & \label{eq:due1}  \\
 & W_m^T AV_m  = I_m & \label{eq:tre1}  \\
 & W_m^T A^2 V_m  = T_m & \label{eq:quattro1}
\end{align}

\end{Proposition}

\textbf{Proof.}~See~\cite{hgzl:09}.

\begin{flushright}
 $\blacksquare$
\end{flushright}

A characterization of Algorithm~\ref{alg:biorth} in terms of projections into
relevant Krylov subspaces is derived in the following 

\begin{Corollary}~\label{cor:biorth}
The biconjugate Lanczos $A$-orthonormalization method ideally builds up a pair 
of biconjugate A-orthonormal bases for the dual Krylov subspaces $K_m(A; v_1)$ 
and $A^T K_m(A^T;w_1)$. The matrix $T_m$ is the projection of 
the matrix $A$ onto the corresponding Krylov subspaces.
\end{Corollary}

\begin{flushright}
 $\blacksquare$
\end{flushright}

\section{The Lanczos $A$-biorthonormalization method for general linear systems}~\label{sec:due}

From the recursions defined in Eqns~(\ref{rec:1})-(\ref{rec:2}) and the characterization 
given by Corollary~(\ref{cor:biorth}), we derive a Lanczos $A$-biorthonormalization method 
for general linear systems along the following lines.

\begin{description}
 \item [STEP~1] Run Algorithm~\ref{alg:biorth} $m$ steps for some
user specified $m \ll n$ and generate Lanczos-type matrices 
$V_m=[v_1,v_2,\ldots,v_m]$, $W_m=[w_1,w_2,\ldots,w_m]$ and the 
tridiagonal projection matrix $T_m$ defined in Proposition~\ref{theo:biorth}.

 \item [STEP~2] Compute the approximate solution $x_m$ that belongs to 
the Krylov subspace $x_0+K_m(A;v_1)$ by imposing the residual be orthogonal
to the constraint subspace $L_m \equiv A^T K_m(A^T;w_1)$

\[
r_m = b - Ax_m~\bot~L_m,
\]

or equivalently in matrix formulation

\begin{equation}\label{eq:proj}
    \left( {A^T W_m } \right)^T \left( {b - Ax_m } \right) = 0.
\end{equation}

Recall that the approximate solution has form 

\begin{equation}\label{eq:basis}
  x_m=x_0+V_m y_m,
\end{equation}

so that by simple substitution and computation with~Eqns~(\ref{eq:quattro1}-\ref{eq:basis}) 
we obtain a tridiagonal system to solve for $y_m$,

\begin{equation}\label{eq:tm}
  T_m y_m  = \beta e_1,~~\beta  = \left\| {r_0 } \right\|_2.
\end{equation}

 \item [STEP~3] Compute the new residual and if convergence is observed, terminate 
the computation. Otherwise, enlarge the Krylov subspace and repeat the process again.

\end{description}

The whole iterative scheme is sketched in Algorithm~\ref{alg:2sided}. 
It solves not only the original linear system $Ax=b$ but, implicitly, also 
a dual linear system $A^T x^* = b^*$ with $A^T$.
Analogously, we can derive the counterpart of Algorithm~\ref{alg:2sided} for the solution of the
corresponding dual system $A^T x^* = b^*$, where the dual approximation $x_m^*$ is sought from the
affine subspace $x_0^*  + K_m \left( {A^T ,w_1 } \right)$ of dimension $m$ satisfying 

\[
b^*  - A^T x_m^*~~\bot~~AK_m \left( {A,v_1 } \right).
\]

Denote $r_0^*  = b^*  - A^T x_0^* $ and $\beta ^*  = \left\| {r_0^* } \right\|_2 $. If 
$w_1  = \frac{{r_0^* }}{{\beta ^* }}$ and $v_1$ is chosen properly such that $\left\langle {v_1 ,Aw_1 } \right\rangle  = 1$,
then the counterparts of~(\ref{eq:proj}-\ref{eq:tm}) are the following

\begin{align}
 & \left( {AV_m } \right)^T \left( {b^*  - A^T x_m^* } \right) = 0 \label{eq:8} \\ 
 & x_m^*  = x_0^*  + W_m y_m^* , \label{eq:9} \\ 
 & T_m^T y_m^*  = \beta ^* e_1 \label{eq:10} 
\end{align}

where $V_m, W_m$ and $T_m$ are defined in Proposition~\ref{theo:biorth} and $y_m^*  \in \R^n $
is the coefficient vector of the dual linear combination. In Section~\ref{sec:tre}
we will derive a formulation that does not require multiplications by $A^T$.

\medskip

\begin{algorithm}
\caption{\label{alg:2sided}{\it Two-sided Biconjugate A-Orthonormalization method}.}
\begin{algorithmic}[1]
  \STATE{Compute $r_0=b-Ax_0$ for some initial guess $x_0$ and set $\beta  = \left\| {r_0 } \right\|_2 $.}
  \STATE{Start with $v_1  = \frac{{r_0 }}{\beta }$ and choose $\omega _1 $ such that $\left\langle {\omega _1 ,Av_1 } \right\rangle  = 1$, (for example, $\omega _1  = \frac{{Av_1 }}{{\left\| {Av_1 } \right\|_2^2 }}$).}
  \STATE{Generate the Lanczos-type vectors 
$V_m  = \left[ {v_1 ,v_2 , \ldots ,v_m } \right]$ and 
$W_m  = \left[ {w_1 ,w_2 , \ldots ,w_m } \right]$ as well as the tridiagonal matrix $T_m$ 
defined in Proposition~\ref{theo:biorth} by running Algorithm~1 $m$ steps.}
  \STATE{Compute $y_m  = T_m^{ - 1} \left( {\beta e_1 } \right)$ and $x_m  = x_0  + V_m y_m $.}
\end{algorithmic}
\end{algorithm}

The approximation $x_m$ and the dual 
approximation $x_m^*$ can be updated respectively from $x_{m-1}$ and $x_{m-1}^*$ at each step. 
Assume the $LU$ decomposition of the tridiagonal matrix $T_m$ is

\[
T_m  = L_m U_m ,
\]

substituting which into~(\ref{eq:basis}-\ref{eq:tm}) and~(\ref{eq:9}-\ref{eq:10}) gives respectively

\begin{align*}
x_m & = x_0  + V_m \left( {L_m U_m } \right)^{ - 1} \left( {\beta e_1 } \right) \\
    & = x_0  + V_m U_m^{ - 1} L_m^{ - 1} \left( {\beta e_1 } \right) \\
    & = x_0  + P_m z_m , \\
x_m^* & = x_0^*  + W_m \left( {U_m^T L_m^T } \right)^{ - 1} \left( {\beta^* e_1 } \right) \\
    & = x_0^*  + W_m \left( {L_m^T} \right)^{ - 1} \left( {U_m^T} \right)^{ - 1} \left( {\beta^* e_1 } \right) \\
    & = x_0^*  + P_m^* z_m^* ,
\end{align*}

where $P_m  = V_m U_m^{ - 1} ,z_m  = L_m^{ - 1} \left( {\beta e_1 } \right)$, and
$P_m^*  = W_m\left( {L_m^T } \right)^{ - 1} ,z_m^*  = \left( {U_m^T } \right)^{ - 1} \left( {\beta ^* e_1 } \right)$.
Using the same argument as in the derivation of DIOM 
from IOM Algorithm in~\cite{saad:96}-Chapter~6, we easily derive the relations

\begin{align*}
& x_m  = x_{m - 1}  + \zeta _m p_m , \\ 
& x_m^*  = x_{m - 1}^*  + \zeta _m^* p_m^* , \\ 
\end{align*}

where $\xi_m$ and $\xi_m^*$ are coefficients, $p_m$ and $p_m^*$ are the corresponding
$m$th column vectors in $P_m$ and $P_m^*$ defined above, termed as the $m$th primary and
dual direction vectors, respectively.

Observe that the pairs of the primary and dual direction vectors form a biconjugate 
$A^2$-orthonormal set, 
\ie,~$\left\langle {p_i^* ,A^2 p_j } \right\rangle  = \delta _{i,j} \left( {1 \le i,j \le m} \right)$,
which follows clearly from 

\begin{align*}
\left( {P_m^* } \right)^T A^2 P_m & = \left( {W_m \left( {L_m^T } \right)^{ - 1} } \right)^T A^2 V_m U_m^{ - 1} \\
 & = L_m^{ - 1} W_m^T A^2 V_m U_m^{ - 1} \\
 & = L_m^{ - 1} T_m U_m^{ - 1} \\
 & = L_m^{ - 1} L_m U_m U_m^{ - 1} \\
 & = I_m
\end{align*}

with~(\ref{eq:quattro1}). 
In addition, the $m$th primary residual vector $r_m=b-Ax_m$ and the $m$th dual residual vector
$r_m^*=b^*-Ax_m^*$ can be expressed as

\begin{align}
& r_m  =  - \delta _{m + 1} e_m^T y_m v_{m + 1} , \label{eq:drsd1} \\ 
& r_m^*  =  - \beta _{m + 1} e_m^T y_m^* w_{m + 1} , \label{eq:drsd2} 
\end{align}

by simple computation 
with~(\ref{eq:uno1}-\ref{eq:due1},\ref{eq:basis}-\ref{eq:tm},\ref{eq:9}-\ref{eq:10}).
Eqns~(\ref{eq:drsd1}) and (\ref{eq:drsd2}) together with~(\ref{eq:tre1}) reveal that 
the primary and dual residual vectors satisfy the biconjugate $A$-orthogonal conditions, 
\ie~$\left\langle {r_i^* ,Ar_j } \right\rangle  = 0$ for $i \ne j$.

\smallskip

It is known that the constraint subspace in an oblique projection method is different 
from the search subspace and may be totally unrelated to it. The distinction is rather 
important and gives rise to different types of algorithms~\cite{saad:96}. The biconjugate 
$A$-orthogonality between the primary and dual residual vectors and the biconjugate 
$A^2$-orthonormality between the primary and dual direction vectors reveal and suggest 
alternative choices for the constraint subspace. This idea helps to devise the BiCOR 
method described in the coming section.

\FloatBarrier

\section{The BiCOR method}

Proceeding in a similar way like the Conjugate Gradient and its variants CG, CR, COCG, 
BiCGCR, COCR and BiCR, given an initial guess $x_0$ to the considered linear system 
$Ax=b$, we derive algorithms governed by the following coupled two-term recurrences 

\begin{align}
 &  r_0  = b - Ax_0 ,~p_0  = r_0 , \label{eq:uno}  \\ 
 &  x_{j + 1}  = x_j  + \alpha _j p_j , \label{eq:due}  \\ 
 &  r_{j + 1}  = r_j  - \alpha _j Ap_j , \label{eq:tre}  \\ 
 &  p_{j + 1}  = r_{j + 1}  + \beta _j p_j ,~\textrm{for}~j = 0,1, \ldots \label{eq:quattro}  
\end{align}
   
where $r_j = b-Ax_j$ is the $j$th residual vector and $p_j$ is the $j$th search 
direction vector. Denoting $L_m$ the underlying constraints subspace, the parameters 
$\alpha_i$,  $\beta_i$ can be determined by imposing the following orthogonality 
conditions:

\begin{equation}\label{eq:17}
r_{j + 1}~\bot~L_m \textrm{~and~} Ap_{j + 1}~\bot~L_m.
\end{equation}

For concerned real unsymmetric linear systems

\begin{equation}
 Ax=b,
\end{equation}

an expanded choice for the constraint subspace is $L_m = A^T K_m (A^T ,r_0^* )$, 
where $r_0^*$ is chosen to be equal to $P(A)r_0$,
with $P(t)$ an arbitrary polynomial of certain degree with respect to the variable
$t$ and $p_0^*=r_0^*$. It should be noted that the optimal choice for the involved
polynomial is in general not easily obtainable and requires some expertise and
artifice. This aspect needs further research. When there is no ambiguity or other
clarification, a specific default choice for $L_m$ with $r_0^*=Ar_0$ is adopted in
the numerical experiments against the other popular choice for $L_m$ with $r_0^*=r_0$,
see e.g.~\cite{sosz:09,soga:06}. It is important to note that the scalars 
~$\alpha_j, \beta_j (j=0,1,\ldots)$ in the recurrences~(\ref{eq:due}-\ref{eq:quattro}) 
are different from those produced by Algorithm~\ref{alg:biorth}. The search direction 
vectors $p_j's$ here are multiples of the primary direction vectors $p_j's$ defined 
in Section~\ref{sec:due}.
The coupled two-term recurrences for the 
$(j+1)$th shadow residual vector $r_{j+1}^*$ and the associated $(j+1)$th shadow search
direction $p_{j+1}^*$ can be augmented by similar relations to~(\ref{eq:tre}-\ref{eq:quattro})
as follows:

\begin{align}
 & r_{j + 1}^*  = r_j^*  - \alpha _j A^T p_j^* , \label{eq:cinque} \\ 
 & p_{j + 1}^*  = r_{j + 1}^*  + \beta _j p_j^* , \textrm{~for~} j = 0,1, \ldots  \label{eq:sei} 
\end{align}

where $ \alpha _j $ and $ \beta _j $ are the conjugate complex of 
$\alpha _j $ and $\beta _j $ in~(\ref{eq:tre}-\ref{eq:quattro}), correspondingly.
Then with a certain polynomial $P(t)$ with respect to $t$,~(\ref{eq:17}) explicitly
reads

\begin{equation}\label{eq:21}
r_{j + 1} ~ \bot ~ A^T K_m \left( {A^T ,r_0^* } \right) \textrm{~and~} 
Ap_{j + 1} ~ \bot ~ A^T K_m \left( {A^T ,r_0^* } \right) \textrm{~with~} 
r_0^*=P(A)r_0
\end{equation}

which can be reinterpreted from a practical point of view as

\begin{itemize}
 \item the residual vectors $r_i$'s and the shadow residual vectors $r_j^*$'s
are biconjugate $A$-orthogonal to each other,~\ie, 
$\left\langle {r_j^* ,Ar_i } \right\rangle  = \left\langle {A^T r_j^* ,r_i } \right\rangle  = 0,$
for $i \ne j$;
 \item the search direction vectors $p_i's$ and the shadow search direction vectors
$p_j^*$'s form an $A^2$-biconjugate set, \ie, 
$\left\langle {p_j^* ,A^2 p_i } \right\rangle  = \left\langle {A^T p_j^* ,Ap_i } \right\rangle  = \left\langle {\left( {A^T } \right)^2p_j^* ,p_i } \right\rangle  = 0,$
for $i \ne j$.
\end{itemize}

These facts are already stated in the latter part of Section~\ref{sec:due}.
Therefore, we possess the conditions to determine the scalars $\alpha_j$
and $\beta_j$ by imposing the corresponding biorthogonality and biconjugacy 
conditions~(\ref{eq:21}) into~(\ref{eq:tre}-\ref{eq:quattro},\ref{eq:cinque}-\ref{eq:sei}). 
We use extensively the algorithmic 
schemes introduced in~\cite{barr:95} for descriptions of the present algorithms.

Making the inner product of $A^T r_j^*$ and $r_{j+1}$ as defined by~(\ref{eq:tre}) yields
\[
\left\langle {A^T r_j^* ,r_{j + 1} } \right\rangle  = \left\langle {A^T r_j^* ,r_j  - \alpha _j Ap_j } \right\rangle  = 0,
\]
with the biconjugate $A$-orthogonality between $r_{j+1}$ and $r_J^*$, further resulting in
\[
\alpha _j  = \frac{{\left\langle {A^T r_j^* ,r_j } \right\rangle }}{{\left\langle {A^T r_j^* ,Ap_j } \right\rangle }},
\]
where the denominator of the above right-hand side can be further modified as
\[
\left\langle {A^T r_j^* ,Ap_j } \right\rangle  = \left\langle {A^T p_j^*  - \beta _{j - 1} A^T p_{j - 1}^* ,Ap_j } \right\rangle  = \left\langle {A^T p_j^* ,Ap_j } \right\rangle ,
\]
because $p_{j-1}^*$ and $p_j$ are $A^2$-biconjugate. Then
\begin{equation}\label{eq:22}
\alpha _j  = \frac{{\left\langle {A^T r_j^* ,r_j } \right\rangle }}{{\left\langle {A^T p_j^* ,Ap_j } \right\rangle }} = \frac{{\left\langle {r_j^* ,Ar_j } \right\rangle }}{{\left\langle {A^T p_j^* ,Ap_j } \right\rangle }} 
\end{equation}

Similarly, writing that $p_{j+1}$ as defined by~(\ref{eq:quattro}) is $A^2$-biconjugate to 
$p_j^*$ yields

\[
\left\langle {p_j^* ,A^2 p_{j + 1} } \right\rangle  = \left\langle {A^T p_j^* ,Ap_{j + 1} } \right\rangle  = \left\langle {A^T p_j^* ,Ar_{j + 1}  + \beta _j Ap_j } \right\rangle  = 0,
\]

giving

\[
\beta _j  =  - \frac{{\left\langle {A^T p_j^* ,Ar_{j + 1} } \right\rangle }}{{\left\langle {A^T p_j^* ,Ap_j } \right\rangle }} =  - \alpha _j \frac{{\left\langle {A^T p_j^* ,Ar_{j + 1} } \right\rangle }}{{\left\langle {r_j^* ,Ar_j } \right\rangle }}
\]

with $\alpha_j$ computed in~(\ref{eq:22}).

Observe from~(\ref{eq:cinque}) that 
\[
 - \alpha _j A^T p_j^*  = r_{j + 1}^*  - r_j^* 
\]
and therefore,

\begin{equation}\label{eq:23}
\beta _j  = \frac{{\left\langle { - \alpha _j A^T p_j^* ,Ar_{j + 1} } \right\rangle }}{{\left\langle {r_j^* ,Ar_j } \right\rangle }} = \frac{{\left\langle {r_{j + 1}^*  - r_j^* ,Ar_{j + 1} } \right\rangle }}{{\left\langle {r_j^* ,Ar_j } \right\rangle }} = \frac{{\left\langle {r_{j + 1}^* ,Ar_{j + 1} } \right\rangle }}{{\left\langle {r_j^* ,Ar_j } \right\rangle }}
\end{equation}

because of the biconjugate $A$-orthogonality of $r_j^*$ and $r_{j+1}$. Putting these
relations~(\ref{eq:uno}-\ref{eq:quattro},\ref{eq:cinque}-\ref{eq:sei},\ref{eq:22}-\ref{eq:23}) 
together and taking the strategy of reducing the number of matrix-vector multiplications 
by introducing an auxiliary vector recurrence and changing variables,
together lead to the BiCOR method. The pseudocode for the left
preconditioned BiCOR method with a preconditioner $M$ is given in the following
Algorithm~\ref{alg:bicor}.

\begin{algorithm}
\caption{\label{alg:bicor}{\it Left preconditioned BiCOR method}.}
\begin{algorithmic}[1]
\STATE Compute $r_0=b-Ax_0$ for some initial guess $x_0$.
\STATE Choose $r_0^*=P(A)r_0$ such that $\left\langle {r_0^* ,Ar_0 } \right\rangle  \ne 0$, where $P(t)$ is a polynomial in $t$. (For example, $r_0^*=Ar_0$).
\FOR{$j=1,2,\ldots$}
  \STATE{\textbf{solve}~$Mz_{j - 1}  = r_{j - 1}$}
  \IF{j=1}
      \STATE{\textbf{solve}~$M^T z_0^*  = r_0^* $}
  \ENDIF
  \STATE $\hat z = Az_{j - 1} $
  \STATE $\rho _{j - 1}  = \left\langle {z_{j - 1}^* ,\hat z} \right\rangle $
  \STATE{\textbf{if}~$\rho _{j - 1}  = 0$,~\textbf{method fails}}
  \IF{$j=1$}
  \STATE $p_0=z_0$
  \STATE $p_0^*=z_0^*$
  \STATE $q_0=\hat z$
  \ELSE
  \STATE $\beta _{j - 2}  = \rho _{j - 1} /\rho _{j - 2} $
  \STATE $p_{j - 1}  = z_{j - 1}  + \beta _{j - 2} p_{j - 2} $
  \STATE $p_{_{j - 1} }^*  = z_{j - 1}^*  + \beta _{j - 2} p_{j - 2}^* $
  \STATE $q_{_{j - 1} }  = \hat z + \beta _{j - 2} q_{j - 2} $
  \ENDIF
  \STATE $q_{_{j - 1} }^*  = A^T p_{_{j - 1} }^* $
  \STATE{\textbf{solve} $M^T u_{j - 1}^*  = q_{j - 1}^*$}
  \STATE $\alpha _{j - 1}  = \rho _{j - 1} /\left\langle {u_{j - 1}^* ,q_{j - 1} } \right\rangle $
  \STATE $x_j  = x_{j - 1}  + \alpha _{j - 1} p_{j - 1} $
  \STATE $r_j  = r_{j - 1}  - \alpha _{j - 1} q_{j - 1} $
  \STATE $z_j^*  = z_{j - 1}^*  - \alpha _{j - 1} u_{j - 1}^* $
  \STATE check convergence; continue if necessary
\ENDFOR
\end{algorithmic}
\end{algorithm}

\FloatBarrier

\section{A transpose-free variant of the BiCOR method: the CORS method}

Exploiting similar ideas to the ingenious derivation of the CGS method~\cite{sonn:89},
one variant of the BiCOR method can be developed which does not require matrix-vector 
products by $A^T$. The new algorithm is referred to as CORS and is derived using a 
different polynomial representation of the residual with the hope of increasing the 
effectiveness of BiCOR in certain circumstances.
First, the CORS method follows exactly a similar way in~\cite{sonn:89} for the 
derivation of the CGS method while taking the strategy of reducing the number
of matrix-vector multiplications by introducing auxiliary vector recurrences
and changing variables. In Algorithm~\ref{alg:bicor}, by simple induction, the 
polynomial representations of the vectors $r_j$, $r_j^*$, $p_j$, $p_j^*$ at step 
$j$ can be expressed as follows

\begin{align*}
  &  r_j  = \phi _j (A)r_0,~~p_j  = \pi _j (A)r_0, \\
  &  r_j^*  = \phi _j (A^T)r_0^*,~~p_j^*  = \pi _j (A^T)r_0^*,	
\end{align*}

where $\phi _j $ and $\pi _j $ are Lanczos-type polynomials of degree less
than or equal to $j$ satisfying $\phi _j (0) = 1$.
Substituting these corresponding polynomial representations 
into~(\ref{eq:22},\ref{eq:23}) gives

\begin{align*}
& \alpha _j  = \frac{{\left\langle {r_j^* ,Ar_j } \right\rangle }}{{\left\langle {A^T p_j^* ,Ap_j } \right\rangle }} = \frac{{\left\langle {\varphi _j \left( {A^T } \right)r_0^* ,A\varphi _j \left( A \right)r_0 } \right\rangle }}{{\left\langle {A^T \pi _j \left( {A^T } \right)r_0^* ,A\pi _j \left( A \right)r_0 } \right\rangle }} = \frac{{\left\langle {r_0^* ,A\varphi _j^2 \left( A \right)r_0 } \right\rangle }}{{\left\langle {r_0^* ,A^2 \pi _j^2 \left( A \right)r_0 } \right\rangle }}, \\
& \beta _j  = \frac{{\left\langle {r_{j + 1}^* ,Ar_{j + 1} } \right\rangle }}{{\left\langle {r_j^* ,Ar_j } \right\rangle }} = \frac{{\left\langle {\varphi _{j + 1} \left( {A^T } \right)r_0^* ,A\varphi _{j + 1} \left( A \right)r_0 } \right\rangle }}{{\left\langle {\varphi _j \left( {A^T } \right)r_0^* ,A\varphi _j \left( A \right)r_0 } \right\rangle }} = \frac{{\left\langle {r_0^* ,A\varphi _{j + 1}^2 \left( A \right)r_0 } \right\rangle }}{{\left\langle {r_0^* ,A\varphi _j^2 \left( A \right)r_0 } \right\rangle }}.
\end{align*}

Also, note from~(\ref{eq:tre},\ref{eq:quattro}) that $\phi _j $ and $\pi _j $ can be 
expressed by the following recurrences

\begin{align*}
& \phi _{j + 1} (t) = \phi _j (t) - \alpha _j t\pi _j (t), \\ 
& \pi _{j + 1} (t) = \phi _{j + 1} (t) + \beta _j \pi _j (t). \\ 
\end{align*}

By some algebraic computation with the help of the induction relations between
$\phi _j $ and $\pi _j $ and the strategy of reducing operations mentioned
above, the desired CORS method can be obtained. The pseudocode for the 
resulting left preconditioned CORS method with a preconditioner $M$ can be
represented by the following scheme.
In many cases, the CORS method is amazingly competitive with the BiCGSTAB method 
(see e.g. Section~\ref{sec:tre}). However, the CORS method, like the CGS, SCGS, 
and the CRS methods, is based on squaring the residual polynomial. In cases of 
irregular convergence, this may lead to substantial build-up of rounding errors 
and worse approximate solutions, or possibly even overflow (see e.g. example). 
For discussions on this effect and its consequences,
see~\cite{sisz:07,saad:96,vdvo:92,barr:95,slva:95}.

\begin{algorithm}
\caption{\label{alg:cors}{\it Left preconditioned CORS method}.}
\begin{algorithmic}[1]
\STATE Compute $r_0=b-Ax_0$ for some initial guess $x_0$.
\STATE Choose $r_0^*=P(A)r_0$ such that $\left\langle {r_0^* ,Ar_0 } \right\rangle  \ne 0$, where $P(t)$ is a polynomial in $t$. (For example, $r_0^*=Ar_0$).
\FOR{$j=1,2,\ldots$}
  \STATE{\textbf{solve} $Mz_{j - 1}  = r_{j - 1} $}
  \STATE $\hat r = Az_{j - 1} $
  \STATE $\rho _{j - 1}  = \left\langle {r_0^* ,\hat r} \right\rangle $
  \STATE{\textbf{if}{$\rho _{j - 1}  = 0$}, \textbf{method fails}}
  \IF{$j  = 1$}
  \STATE{$e_0  = r_0 $}
  \STATE{solve $Mze_0  = e_0 $}
  \STATE{$d_0  = \hat r$}
  \STATE{$q_0  = \hat r$}
  \ELSE
  \STATE{$\beta _{j - 2}  = \rho _{j - 1} /\rho _{j - 2} $}
  \STATE{$e_{j - 1}  = r_{j - 1}  + \beta _{j - 2} h_{j - 2} $}
  \STATE{$ze_{j - 1}  = z_{j - 1}  + \beta _{j - 2} zh_{j - 2} $}
  \STATE{$d_{j - 1}  = \hat r + \beta _{j - 2} f_{j - 2} $}
  \STATE{$q_{j - 1}  = d_{j - 1}  + \beta _{j - 2} \left( {f_{j - 2}  + \beta _{j - 2} q_{j - 2} } \right)$}
  \ENDIF
  \STATE{\textbf{solve}~$Mq = q_{j - 1} $}
  \STATE{$\hat q = Aq$}
  \STATE{$\alpha _{j - 1}  = \rho _{j - 1} /\left\langle {r_0^* ,\hat q} \right\rangle $}
  \STATE{$h_{j - 1}  = e_{j - 1}  - \alpha _{j - 1} q_{j - 1} $}
  \STATE{$zh_{j - 1}  = ze_{j - 1}  - \alpha _{j - 1} q$}
  \STATE{$f_{j - 1}  = d_{j - 1}  - \alpha _{j - 1} \hat q$}
  \STATE{$x_j  = x_{j - 1}  + \alpha _{j - 1} \left( {2ze_{j - 1}  - \alpha _{j - 1} q} \right)$}
  \STATE{$r_j  = r_{j - 1}  - \alpha _{j - 1} \left( {2d_{j - 1}  - \alpha _{j - 1} \hat q} \right)$}
  \STATE check convergence; continue if necessary
\ENDFOR
\end{algorithmic}
\end{algorithm}

\FloatBarrier

\section{Numerical experiments}~\label{sec:tre}

We initially illustrate the numerical behavior of the proposed algorithms on a set of 
sparse linear systems. Although the focus of the paper is on real unsymmetric
problems, we also report on experiments on complex systems. The Lanczos biconjugate
$A$-orthonormalization procedure is straightforward to generalize to complex matrices,
see e.g.~\cite{hgzl:09}. 
We select matrix problems of different size, from small (a few tens thousand
unknowns) to large (more than one million unknowns), and arising from various application 
areas. The test problems are extracted from the University of Florida~\cite{davi:94} 
matrix collection, except the two SOMMEL problems that are made available at Delft University.
The characteristics of the model problems are illustrated in Tables~(\ref{tab:real},\ref{tab:complex}), 
and the numerical results are shown in Tables~\ref{tab:noshadow}. 
The experiments are carried out using double precision floating point arithmetic 
in MATLAB 7.7.0 on a PC equipped with a Intel(R) Core(TM)2 Duo CPU P8700 running 
at 2.53GHz, and with 4 GB of RAM. 
We report on number of iterations (referred to as $Iters$), CPU consuming time 
in seconds (referred to as $CPU$), $\log_{10}$ of the final true relative 
residual 2-norms defined
$\log_{10} (||b - Ax_n ||_2 / ||r_0 ||_2)$ (referred to as 
$TRR$). 
The stopping criterium used here is that the 2-norm of the residual be reduced by a factor
$TOL$ of the 2-norm of the initial residual, i.e., 
$||r_n ||_2 / ||r_0 ||_2 < TOL$, or when $Iters$ exceeded the maximal number of matrix-by-vector
products $MAXMV$. We take $MAXMV$ very large ($10000)$ and $TOL = 10^{-8}$. 
All these tests are started with a zero initial guess. 
The physical right-hand side is not always available for all problems. Therefore
we set $b = Ae$, where $e$ is the $n \times 1$ vector whose elements are all equal to unity, 
such that $x = (1, 1, · · · , 1)^T$ is the exact solution. 
It is stressed that the specific default choice for the constraint space $L_m$ with 
$r_0^* = A r_0$ is adopted in the implementation of the Lanczos biconjugate A-orthonormalization
methods. Finally, a symbol '-' is used to indicate that the method did not
meet the required TOL before $MAXMV$.

We observe that the proposed algorithms enable to solve fairly large problems in a moderate number 
of iterations. The iteration count is always much smaller than the problem dimension, see 
e.g. experiments on the large ATMOSMODJ, ATMOSMODL, KIM2 problems. 
The final approximate solution is generally very accurate.
In all our experiments we did not observe breakdowns in the Lanczos $A$-biorthonormalization
procedure that proves to be remarkably robust.
The setup of the methods does not require parameters; however, there is some freedom in 
the selection of the initial shadow residual $r_0^*$. 
We analyse in Table~\ref{tab:shadow} the effect of a different choice
for the initial shadow residual; performance may slightly change but
the convergence trend is preserved and in general it is not possible to 
predict a priori the effect.
The results indicate that
CORS is significantly  more robust and faster than BiCOR.
The proposed family of solvers exhibits fast convergence, is parameter-free,
is extremely cheap in memory as it is derived from short-term vector recurrences
and does not suffer from the restriction to require a symmetric preconditioner
when it is applied to symmetric systems.
Therefore, it can be a suitable computational tool to use in applications.

In Tables~(\ref{tab:1}-\ref{tab:14}) we illustrate some comparative figures with other popular Krylov solvers,
that are developed on either the Arnoldi and the Lanczos biconjugation methods.
The only intention of these experiments is to show the level of competiteveness
of the presented Lanczos $A$-biorthonormalization method with other popular approaches
for linear systems. We consider GMRES(50), Bi-CG,QMR, CGS, 
Bi-CGSTAB, IDR(4), BiCOR, CORS, BiCR. Indeed some of these algorithms 
depend on parameters and we did not search the optimal value on our problems.
We set the value of restart of GMRES equal to $50$. The memory request for GMRES 
is the matrix+$($restart$+5)n$, so that it remains much larger than the memory
needed for CORS and BiCOR (see Table~\ref{tab:complexity}).
The parameter value for IDR is selected to $4$, that is often used in experiments.
The CORS method proves very fast in comparison, thanks to its low algorithmic 
cost. We observe the robustness of the two algorithms on the STOMMEL1 problem
from Ocean modelling. On this problem, CORS, BICOR and IDR are equally efficient, 
IDR being slightly faster; however, CORS and BICOR are remarkably accurate.

\begin{table}[!ht]
 \begin{center}
  \begin{tabular}{|l|r|l|l|}
    \hline
      Matrix problem & Size & Field & Characteristics  \\
    \hline
     WATER\_TANK & 60,740	     & 3D fluid flow		     & real unsymmetric         \\
     XENON2 &  157,464               & Materials 		      & real unsymmetric        \\
     STOMACH &  213,360              & Electrophysiology              & real unsymmetric       	\\
     TORSO3 &  259,156               & Electrophysiology             & real unsymmetric       	\\
     LANGUAGE &  399,130             & Natural language processing   &  real unsymmetric      	\\
     MAJORBASIS & 160,000            & Optimization                  & real unsymmetric       	\\
     ATMOSMODJ & 1,270,432           & Atmospheric modelling          &  real unsymmetric       \\
     ATMOSMODL & 1,489,752     	     & Atmospheric modelling          &  real unsymmetric       \\
    \hline
  \end{tabular}
  \caption{\label{tab:real}Set of test real matrix problems.}
 \end{center}
\end{table}

\begin{table}[!ht]
 \begin{center}
  \begin{tabular}{|l|r|l|l|}
    \hline
      Matrix problem & Size & Field & Characteristics  \\
    \hline
     M4D2\_unscal & 10,000           & Quantum Mechanics	     &  complex unsymmetric     \\
     WAVEGUIDE3D & 21,036            & Electromagnetics              &  complex unsymmetric   	\\
     VFEM &  93,476                  & Electromagnetics              &  complex unsymmetric     \\
     KIM2 &  456,976                 & Complex mesh                  &  complex unsymmetric   	\\
    \hline
  \end{tabular}
  \caption{\label{tab:complex}Set of test complex matrix problems.}
 \end{center}
\end{table}

\begin{table}
\centering
\begin{tabular}{|c|c|c|}
\hline
${\raise0.7ex\hbox{Solver} \!\mathord{\left/
 {\vphantom {{Solver} {Example}}}\right.\kern-\nulldelimiterspace}
\!\lower0.7ex\hbox{Example}}$ & BiCOR & CORS  \\
\hline
\hline
WATER\_TANK 	& 1711 / 38.66 / -8.01 & 1290 / 28.6 / -8.1001   \\
\hline 
XENON2    	& 1247 / 64.79 / -8.0366 & 711 / 37.87 / -8.1368  \\
\hline
STOMACH       	& 82 / 4.33 / -8.1891 & 40 / 2.36 / -8.1321  \\
\hline
TORSO3       	& 97 / 6.92 / -8.12 & 56 / 4.25 /-8.3181  \\
\hline
LANGUAGE 	& 34 / 2.54 / -9.4099 & 24 / 2.11 / -8.3046  \\
\hline
MAJORBASIS    	& 60 / 2.22 / -8.0077 & 35 / 1.45 / -8.688  \\
\hline
ATMOSMODJ     	& 416 / 120 / -8.1018 &  286 / 110.52 / -8.0685  \\
\hline
ATMOSMODL      	& 273 / 106.36 / -8.2345 & 192 / 87.69 / -8.0312  \\
\hline
M4D2\_unscal	& 2324 / 8.38 / -8.0866 & - / 21.43 / 9.1144  \\
\hline
WAVEGUIDE3D     & 3874 / 38.96 / -8.0508 & 2988 / 38.41 / -8.0485  \\
\hline
VFEM           	& 4022 / 226.26 / -8.0156 & 0.5 / 3.17 / 0.72363  \\
\hline
KIM2      	& 189 / 53.86 / -8.4318 & 105 / 36.99 / -8.0199  \\
\hline
\end{tabular}
\caption{\label{tab:noshadow}Number of iterations, CPU time and $\log_{10}$ of the true residual
on all test examples solved with $TOL = 10^{-8}$ using $r_o^*=A r_0$ for the shadow residual.}
\end{table}

\begin{table}
\centering
\begin{tabular}{|c|c|c|}
\hline
${\raise0.7ex\hbox{Solver} \!\mathord{\left/
 {\vphantom {{Solver} {Example}}}\right.\kern-\nulldelimiterspace}
\!\lower0.7ex\hbox{Example}}$ & BiCOR & CORS  \\
\hline
\hline
WATER\_TANK 	& 1423 / 34.35 / -8.0208   &   1246 / 30.36 / -8.1874    \\
\hline 
XENON2    	& 1164 / 64.02 / -8.0014   & 	730 / 39.25 / -8.1178	\\
\hline
STOMACH         & 81 / 4.37 / -8.1078      & 	50 / 2.96 / -8.0726	\\
\hline
TORSO3       	& 87 / 6.19 / -8.2355      & 	56 / 4.57 / -8.1768	 \\
\hline
LANGUAGE 	& 31 / 2.41 / -8.8309      & 	23 / 2.09 / -8.2032	 \\
\hline
MAJORBASIS    	& 56 / 2.36 / -8.0939	   & 	35 / 1.43 / -8.7153	 \\
\hline
ATMOSMODJ     	& 393 / 99.44 / -8.0175    & 	305 / 92.53 / -8.1162	\\
\hline
ATMOSMODL      	& 279 / 83.86 / -8.0279    & 	213 / 75.95 / -8.1906	\\
\hline
M4D2\_unscal 	& 2370 / 8.58 / -8.0888    & 	5000 / 21.15 / 11.8045	\\
\hline
WAVEGUIDE3D     & 3691 / 37.49 / -8.0015   & 	3009 / 35.66 / -8.0041	\\
\hline
VFEM           	& 3902 / 171.04 / -8.0013  & 	4012 / 232.72 / -8.1491	\\
\hline
KIM2      	& 189 / 51.9 / -8.3155     & 	105 / 36.89 / -8.0861	\\
\hline
\end{tabular}
\caption{\label{tab:shadow}Number of iterations, CPU time and $\log_{10}$ of the true residual
on all test examples solved with $TOL = 10^{-8}$ using $r_o^*=r_0$ for the shadow residual.}
\end{table}

\begin{table}
\begin{minipage}[b]{0.5\linewidth} 
\centering
\begin{tabular}{|c|c|c|c|}
\hline
Method & Iters & CPU & TRR \\
\hline
\hline
GMRES(50) & 9 & 0.39 & -2.4224  \\
\hline 
Bi-CG     & 2145 & 72.17 & -8.0419 \\
\hline
QMR       & 2074 & 78.56 & -8.0006 \\
\hline
CGS       & 2162 & 65.8 & -8.1078 \\
\hline
Bi-CGSTAB & 1232 & 50.52 & -8.0027 \\
\hline
IDR(4)    & 4682 & 66.11 & -8.0242 \\
\hline
BiCOR     & 1711 & 38.66 & -8.01 \\
\hline
CORS      & 1290 & \textbf{28.6} & -8.1001 \\
\hline
BiCR      & 1423 & 38.09 & -8.0208 \\
\hline
\end{tabular}
\caption{\label{tab:1}Comparison results of WATER\_TANK with $TOL = 10^{-8}$.}
\end{minipage}
\hspace{0.5cm} 
\begin{minipage}[b]{0.5\linewidth}
\centering
\begin{tabular}{|c|c|c|c|}
\hline
Method & Iters & CPU & TRR \\
\hline
\hline
GMRES(50) & 1 & 0.3 & -0.91923  \\
\hline
Bi-CG     & 1329 & 104.32 & -8.0035  \\
\hline 
QMR       & 1165 & 104.02 & -8.0014  \\
\hline 
CGS       & 968 & 69.52 & -8.258  \\
\hline 
Bi-CGSTAB & 866.5 & 85.08 & -8.1773  \\
\hline 
IDR(4)    & 2148 & 76.41 & -8.1822  \\
\hline 
BiCOR     & 1247 & 64.79 & -8.0366  \\
\hline 
CORS      & 711 & \textbf{37.87} & -8.1368  \\
\hline 
BiCR      & 1164 & 74.48 & -8.0014  \\
\hline
\end{tabular}
\caption{\label{tab:2}Comparison results of XENON2 with $TOL = 10^{-8}$.}
\end{minipage}
\end{table}

\begin{table}
\begin{minipage}[b]{0.5\linewidth} 
\centering
\begin{tabular}{|c|c|c|c|}
\hline
Method & Iters & CPU & TRR \\
\hline
\hline
GMRES(50) & 75 & 20.31 & -8.4244 \\ 
\hline 
Bi-CG     & 80 & 6.52 & -8.1185 \\
\hline 
QMR       & 82 & 7.75 & -8.2002 \\
\hline 
CGS       & 50 & 3.69 & -8.0306 \\
\hline 
Bi-CGSTAB & 149.5 & 16.1 & -8.0366 \\ 
\hline 
IDR(4)    & 107 & 5.38 & -8.3001 \\
\hline 
BiCOR     & 82 & 4.33 & -8.1891 \\
\hline 
CORS      & 40 & \textbf{2.36} & -8.1321 \\
\hline 
BiCR      & 81 & 5.61 & -8.1078 \\
\hline
\end{tabular}
\caption{\label{tab:3}Comparison results of STOMACH with $TOL = 10^{-8}$.}
\end{minipage}
\hspace{0.5cm} 
\begin{minipage}[b]{0.5\linewidth}
\centering
\begin{tabular}{|c|c|c|c|}
\hline
Method & Iters & CPU & TRR \\
\hline
\hline
GMRES(50) & 87 & 29.94 & -8.2677 \\
 \hline
Bi-CG     & 84 & 8.83 & -8.0773 \\
 \hline
QMR       & 91 & 11.12 & -8.0037 \\
 \hline
CGS       & 78 & 7.55 & -8.1121 \\
 \hline
Bi-CGSTAB & 105.5 & 14.08 & -8.0388 \\
 \hline
IDR(4)    & 113 & 6.66 & -8.0402 \\
 \hline
BiCOR     & 97 & 6.92 & -8.12 \\
 \hline
CORS      & 56 & \textbf{4.25} & -8.3181 \\
 \hline
BiCR      & 87 & 7.82 & -8.2355 \\
 \hline
\end{tabular}
\caption{\label{tab:4}Comparison results of TORSO3 with $TOL = 10^{-8}$.}
\end{minipage}
\end{table}

\begin{table}
\begin{minipage}[b]{0.5\linewidth} 
\centering
\begin{tabular}{|c|c|c|c|}
\hline
Method & Iters & CPU & TRR \\
\hline
\hline
GMRES(50) & 29 & 8.91 & -8.112 \\
\hline 
Bi-CG     & 30 & 3.43 & -8.1769 \\
\hline 
QMR       & 30 & 4.24 & -8.1875 \\
\hline
CGS       & 25 & 2.66 & -10.6201 \\
\hline 
Bi-CGSTAB & 24 & 3.58 & -8.0851 \\
\hline 
IDR(4)    & 38 & 2.76 & -8.8881 \\
\hline 
BiCOR     & 34 & 2.54 & -9.4099 \\
\hline 
CORS      & 24 & \textbf{2.11} & -8.3046 \\
\hline 
BiCR      & 31 & 3.3 & -8.8309 \\
\hline
\end{tabular}
\caption{\label{tab:5}Comparison results of LANGUAGE with $TOL = 10^{-8}$.}
\end{minipage}
\hspace{0.5cm} 
\begin{minipage}[b]{0.5\linewidth}
\centering
\begin{tabular}{|c|c|c|c|}
\hline
Method & Iters & CPU & TRR \\
\hline
\hline
GMRES(50) & 41 & 6.61 & -6.9214 \\
\hline 
Bi-CG     & 59 & 3.3 & -8.2861 \\
\hline 
QMR       & 59 & 3.99 & -8.1854 \\
\hline 
CGS       & 34 & 1.77 & -8.2714 \\
\hline 
Bi-CGSTAB & 28.5 & 2.06 & -8.0297 \\
\hline 
IDR(4)    & 58 & 2.08 & -8.2353 \\
\hline 
BiCOR     & 60 & 2.22 & -8.0077 \\
\hline 
CORS      & 35 & \textbf{1.45} & -8.688 \\
\hline 
BiCR      & 56 & 2.73 & -8.0939 \\
\hline
\end{tabular}
\caption{\label{tab:6}Comparison results of MAJORBASIS with $TOL = 10^{-8}$.}
\end{minipage}
\end{table}

\begin{table}
\begin{minipage}[b]{0.5\linewidth} 
\centering
\begin{tabular}{|c|c|c|c|}
\hline
Method & Iters & CPU & TRR \\
\hline
\hline
GMRES(50) & 175 & 336.46 & -3.02 \\
\hline 
Bi-CG     & 432 & 185.34 & -8.0037 \\
\hline 
QMR       & 432 & 210.6 & -8.0712 \\
\hline 
CGS       & 294 & 132.08 & -8.4498 \\
\hline 
Bi-CGSTAB & 266 & 143.87 & -8.4199 \\
\hline 
IDR(4)    & 540 & 135.58 & -8.2371 \\
\hline 
BiCOR     & 416 & 120 & -8.1018 \\
\hline 
CORS      & 286 & \textbf{110.52} & -8.0685 \\
\hline 
BiCR      & 393 & 166.91 & -8.0175 \\
\hline
\end{tabular}
\caption{\label{tab:7}Comparison results of ATMOSMODJ with $TOL = 10^{-8}$.}
\end{minipage}
\hspace{0.5cm} 
\begin{minipage}[b]{0.5\linewidth}
\begin{tabular}{|c|c|c|c|}
\hline
Method & Iters & CPU & TRR \\
\hline
\hline
GMRES(50) & 71 & 123.98 & -2.5268 \\
\hline 
Bi-CG     & 296 & 141 & -8.4055 \\
\hline 
QMR       & 284 & 168.65 & -8.0294 \\
\hline 
CGS       & 240 & 106.3 & -8.0073 \\
\hline 
Bi-CGSTAB & 170 & 105.96 & -8.0032 \\
\hline 
IDR(4)    & 298 & 114.27 & -8.1764 \\
\hline 
BiCOR     & 273 & 106.36 & -8.2345 \\
\hline 
CORS      & 192 & \textbf{87.69} & -8.0312 \\
\hline 
BiCR      & 279 & 120.81 & -8.0279 \\
\hline
\end{tabular}
\caption{\label{tab:8}Comparison results of ATMOSMODL with $TOL = 10^{-8}$.}
\end{minipage}
\end{table}

\begin{table}
\begin{minipage}[b]{0.5\linewidth} 
\centering
\begin{tabular}{|c|c|c|c|}
\hline
Method & Iters & CPU & TRR \\
\hline
\hline
GMRES(50) & 8850 & 309.41 & -0.86138 \\
\hline 
Bi-CG     & 3742 & 43.09 & -8.0516 \\
\hline 
QMR       & 3479 & 52.07 & -8.1005 \\
\hline 
CGS       & 1286 & 54.39 & -0.53875 \\
\hline 
Bi-CGSTAB & 2300 & 33.26 & -8.2957 \\
\hline 
IDR(4)    & 2889 & \textbf{17.86} & -6.4891 \\
\hline 
BiCOR     & 2397 & 18.43 & -8.019 \\
\hline 
CORS      & 2155 & 17.94 & -7.7396 \\
\hline 
BiCR      & 2583 & 28.33 & -8.6085 \\
\hline
\end{tabular}
\caption{\label{tab:9}Comparison results of STOMMEL1 with $TOL = 10^{-8}$.}
\end{minipage}
\hspace{0.5cm} 
\begin{minipage}[b]{0.5\linewidth}
\begin{tabular}{|c|c|c|c|}
\hline
Method & Iters & CPU & TRR \\
\hline
\hline
GMRES(50) & - & 78.64 & -2.4174 \\
\hline 
Bi-CG     & 1515 & 4.14 & -8.0953 \\
\hline 
QMR       & 1492 & 5.54 & -8.0046 \\
\hline 
CGS       & 2520 & 11.63 & -0.79492 \\
\hline 
Bi-CGSTAB & 899.5 & 3.2 & -8.4985 \\
\hline 
IDR(4)    & 1119 & \textbf{1.61} & -8.3146 \\
\hline 
BiCOR     & 1302 & 2.38 & -8.0514 \\
\hline 
CORS      & 936 & 1.84 & -8.1765 \\
\hline 
BiCR      & 1223 & 3.37 & -8.084 \\
\hline
\end{tabular}
\caption{\label{tab:10}Comparison results of STOMMEL2 with $TOL = 10^{-8}$.}
\end{minipage}
\end{table}

\begin{table}
\begin{minipage}[b]{0.5\linewidth} 
\centering
\begin{tabular}{|c|c|c|c|}
\hline
Method & Iters & CPU & TRR \\
\hline
\hline
GMRES(50) & 3127 & 44.69 & -8.0001 \\
\hline 
Bi-CG     & 2488 & 14.17 & -8.0319 \\
\hline 
QMR       & 2497 & 18.25 & -8.0226 \\
\hline 
CGS       & - & 28.74 & -0.10264 \\
\hline 
Bi-CGSTAB & - & 39.46 & -4.7277 \\
\hline 
IDR(4)    & - & 12.14 & -8.0658 \\
\hline 
BiCOR     & 2324 & \textbf{8.38} & -8.0866 \\
\hline 
CORS      & - & 21.43 & 9.1144 \\
\hline 
BiCR      & 2370 & 12.36 & -8.0888 \\
\hline
\end{tabular}
\caption{\label{tab:11}Comparison results of M4D2\_unscal with $TOL = 10^{-8}$.}
\end{minipage}
\hspace{0.5cm} 
\begin{minipage}[b]{0.5\linewidth}
\begin{tabular}{|c|c|c|c|}
\hline
Method & Iters & CPU & TRR \\
\hline
\hline
GMRES(50) & - & 299.8 & -6.9859 \\
\hline 
Bi-CG     & 3851 & 57.93 & -8.1015 \\
\hline 
QMR       & 3752 & 68.24 & -8.0482 \\
\hline 
CGS       & 3013 & 44.88 & -8.0216 \\
\hline 
Bi-CGSTAB & 4957 & 101.45 & -4.9052 \\
\hline 
IDR(4)    & - & 70.27 & -6.2982 \\
\hline 
BiCOR     & 3874 & 38.96 & -8.0508 \\
\hline 
CORS      & 2988 & \textbf{38.41} & -8.0485 \\
\hline 
BiCR      & 3691 & 49.07 & -8.0015 \\
\hline 
\end{tabular}
\caption{\label{tab:12}Comparison results of WAVEGUIDE3D with $TOL = 10^{-8}$.}
\end{minipage}
\end{table}

\begin{table}
\begin{minipage}[b]{0.5\linewidth} 
\centering
\begin{tabular}{|c|c|c|c|}
\hline
Method & Iters & CPU & TRR \\
\hline
\hline
GMRES(50) & - & 2363.42 & -6.3922 \\
\hline 
Bi-CG     & 4744 & 337.54 & -7.9757 \\
\hline 
QMR       & 3786 & 307.21 & -8.0008 \\
\hline 
CGS       & 4709 & 340.64 & -5.4999 \\
\hline 
Bi-CGSTAB & 1941 & 443.64 & -6.809 \\
\hline 
IDR(4)    & - & 360.89 & 5.9272 \\
\hline 
BiCOR     & 3593 & \textbf{154.87} & -8.0113 \\
\hline 
CORS      & 4022 & 226.26 & -8.0156 \\
\hline 
BiCR      & 3902 & 225.49 & -8.0013 \\
\hline 
\end{tabular}
\caption{\label{tab:13}Comparison results of VFEM with $TOL = 10^{-8}$.}
\end{minipage}
\hspace{0.5cm} 
\begin{minipage}[b]{0.5\linewidth}
\begin{tabular}{|c|c|c|c|}
\hline
Method & Iters & CPU & TRR \\
\hline
\hline
GMRES(50) & 40 & 53.08 & -3.0378 \\
\hline 
Bi-CG     & 189 & 79.99 & -8.0716 \\
\hline 
QMR       & 195 & 94.95 & -8.5199 \\
\hline 
CGS       & 105 & 43.42 & -8.1447 \\
\hline 
Bi-CGSTAB & 133 & 75.68 & -8.0007 \\
\hline 
IDR(4)    & 187 & 46.44 & -8.067 \\
\hline 
BiCOR     & 189 & 53.86 & -8.4318 \\
\hline 
CORS      & 105 & \textbf{36.99} & -8.0199 \\
\hline 
BiCR      & 189 & 66.64 & -8.3155 \\
\hline
\end{tabular}
\caption{\label{tab:14}Comparison results of KIM2 with $TOL = 10^{-8}$.}
\end{minipage}
\end{table}

\begin{table}
  \centering
  \begin{small}
\begin{tabular}{l|ccccc}

    Solver & Type & Products by $A$ &  Products by $A^T$ & Scalar products & Memory   \\
\hline

     BiCOR    & general   & 1    & 1 & 2 & matrix+$10n$  \\
      CORS    & "         & 2    & - & 2 & matrix+$14n$   \\

\end{tabular}
  \caption{Algorithmic cost and memory expenses per iteration for Lanczos $A$-biorthonormalization
methods for linear systems.}\label{tab:complexity}
  \end{small}
\end{table}

\subsection{dense problems}

In the last decade, iterative methods have become widespread also in dense
matrix computation partly due to the development of efficient boundary element 
techniques for engineering and scientific applications.
Boundary integral equations, \ie~integral equations defined on the boundary of the
domain of interest, is one of the largest source of truly dense 
linear systems in computational science~\cite{edel:91,edel:93}. 
Recent efforts to implement efficiently these techniques on massively 
parallel computers have resulted in competitive application codes 
provably scalable to several million discretization points,
e.g. the FISC code developed at University of 
Illinois~\cite{solc:97,soch:98,sofi:98a}, the INRIA/EADS integral equation code 
AS\_ELFIP~\cite{sylv:02,sylv:03}, the Bilkent University code~\cite{ergu:07,ergu:08}, 
urging the quest of robust iterative algorithms in this area. 
Integral formulations of surface scattering and hybrid surface/volume formulations 
yield non-Hermitian linear systems that cannot be solved using the Conjugate Gradient 
(CG) algorithm. The GMRES method~\cite{sasc:86} is broadly employed 
due to its robustness and smooth convergence. 
In particular, its flexible variant FGMRES~\cite{saad:93} has become also very popular 
in combination with inner-outer solution schemes, see e.g. experiments reported 
in~\cite{cdgs:05,maeg:07} for solving very large boundary element equations. 
Iterative methods based on the Lanczos biconjugation method, e.g. BiCGSTAB~\cite{vdvo:92} and 
QMR~\cite{frna:91} are also considered in many studies but they may need many 
more iterations to converge especially on realistic geometries~\cite{magu:07,cdgs:05}. 

In this study we  report on results of experiments with Lanczos A-biorthonormalization 
methods on four dense matrix problems arising from boundary integral equations in 
electromagnetic scattering from large structures. 
An accurate solution of scattering problems is a critical concern in civil aviation 
simulations and stealth technology industry, in the analysis of electromagnetic compatibility,
in medical imaging, and other applications.
The selected linear systems arise from reformulating the Maxwell's equations in the frequency 
domain as the following variational problem:

\medskip

Find the distribution of the surface current $\vec{j}$ induced by an incoming radiation,
such that for all tangential test functions $\vec{j}^t$ we have

\begin{multline}\label{eq:EFIE}
\int_\Gamma  {\int_\Gamma  {G(|y - x|)} } \left( {\vec j(x) \cdot
\vec j^t (y) - \frac{1}{{k^2 }} {\rm div}_\Gamma  \vec j(x) \cdot
{\rm div}_\Gamma  \vec j^t (y)} \right)dxdy =   \\
 = \frac{i}{{kZ_0 }}\int_\Gamma {\vec E_{inc} (x)}  \cdot \vec j^t
(x)dx.
\end{multline}

\begin{figure}[!ht]
 \centering
 \subfigure[Sphere]
 { \label{fig:ex5}
   \includegraphics[width=0.45\textwidth]{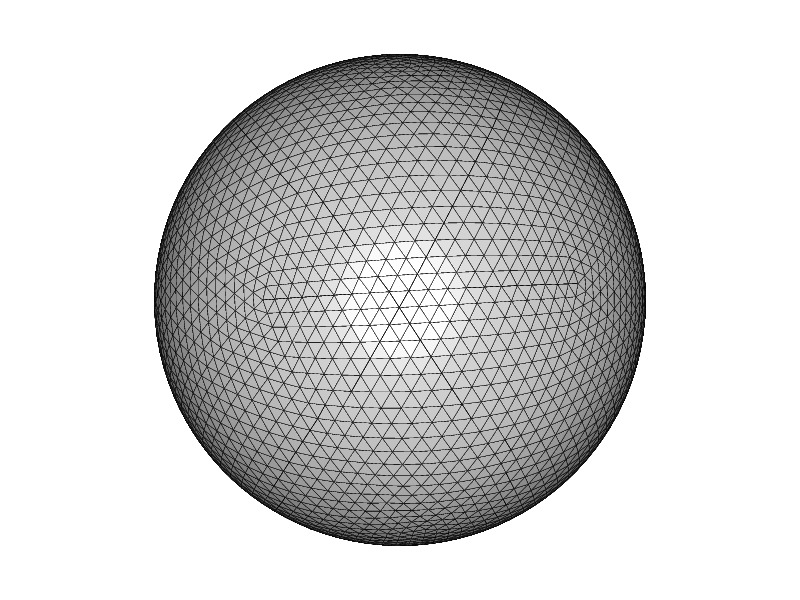}
 }
 \subfigure[Satellite]
 { \label{fig:ex6}
   \includegraphics[width=0.45\textwidth]{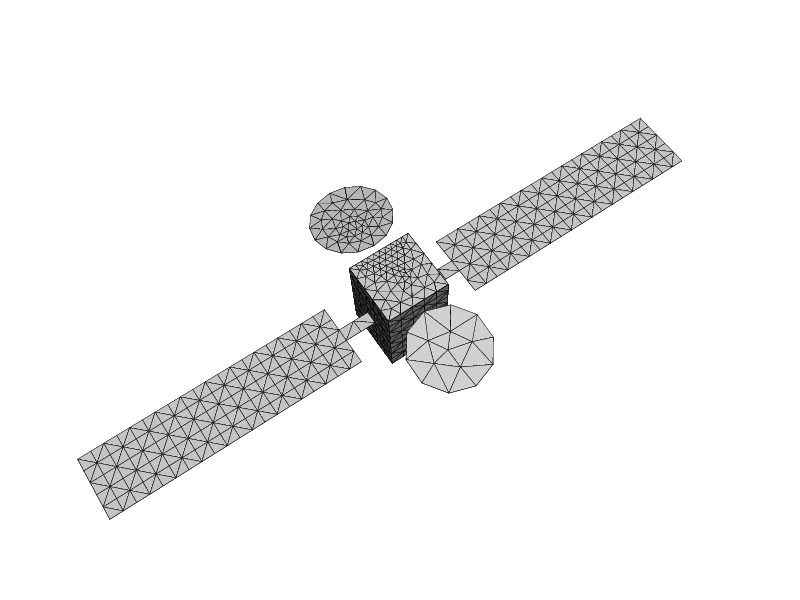}
 }\\
 \subfigure[Jet prototype]
 { \label{fig:ex7}
   \includegraphics[width=0.45\textwidth]{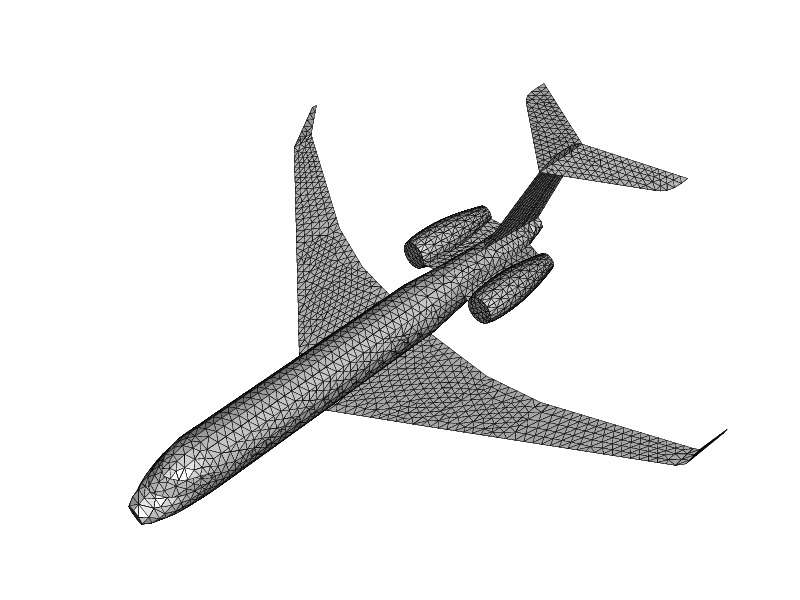}
 }
 \subfigure[Airbus A318 prototype]
 { \label{fig:ex8}
   \includegraphics[width=0.45\textwidth]{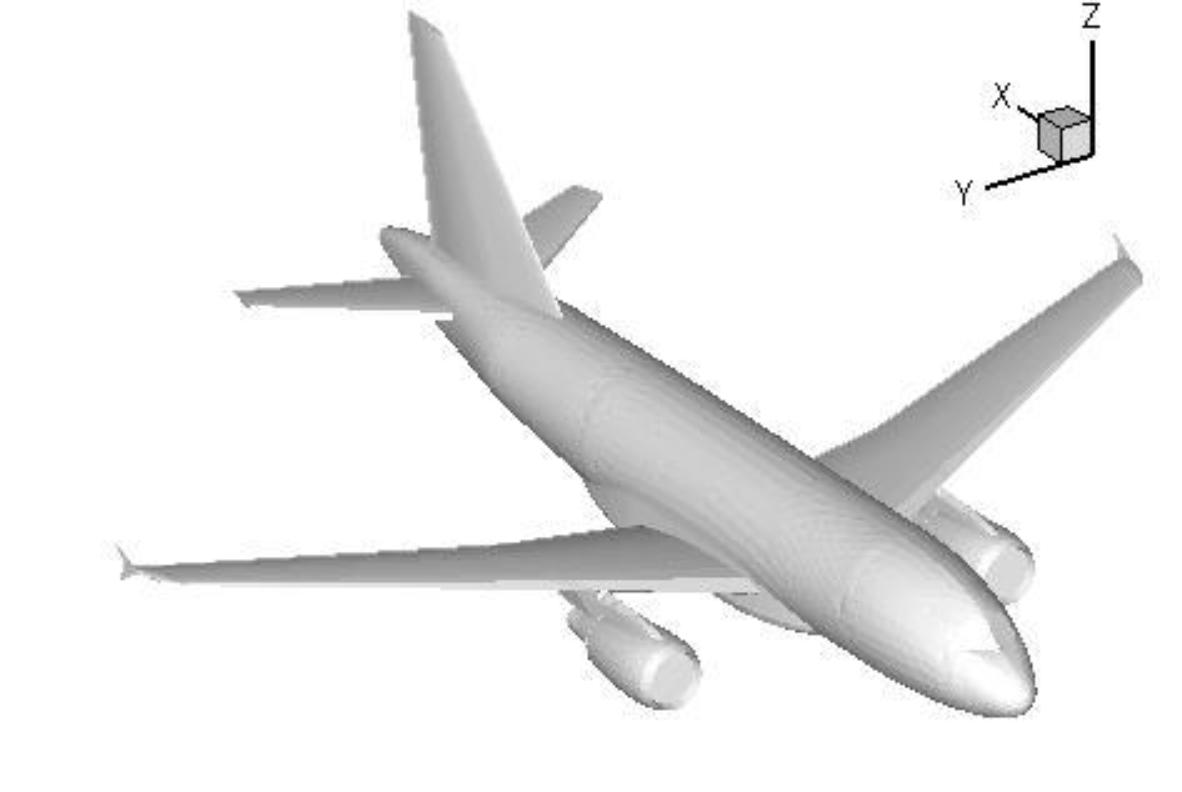}
 }\\
 \caption{\label{fig:mesh} Meshes associated with test examples.}
\end{figure}

We denote by $G(|y-x|)=\displaystyle \frac{e^{ik|y-x|}}{4 \pi |y-x|}$ the Green's function of 
Helmholtz equation, $\Gamma$ the boundary of the object, $k$ the wave number and 
$Z_0  = \sqrt {\mu _0 /\varepsilon _0 } $ the characteristic impedance of vacuum ($\epsilon$ is 
the electric permittivity and $\mu$ the magnetic permeability).
Boundary element discretizations of~(\ref{eq:EFIE}) over a mesh containing $n$ edges produce
dense complex non-Hermitian systems $Ax=b$. The set of unknowns are associated with the 
vectorial flux across an edge in the mesh, while the right-hand side varies with the 
frequency and the direction of the illuminating wave.
When the number of unknowns $n$ is related to the wavenumber, the iteration count of iterative solvers 
typically increase as $\op(n^{0.5})$~\cite{soch:98}. 
Eqn~(\ref{eq:EFIE}) is known as the Electric Field Integral Equation.
Other possible integral models, e.g. the Combined Field Integral
Equation or the Magnetic Field Integral Equation, only apply to closed surfaces 
and are easier to solve by iterative methods~\cite{magu:07}. 
Therefore, we stick with Eqn.~(\ref{eq:EFIE}).
We report the characteristics of the linear systems in Table~\ref{tab:problems} 
and we depict the corresponding geometries in Figure~\ref{fig:mesh}.
In addition to BiCOR and CORS, we consider complex versions of BiCGSTAB, QMR, 
GMRES. Indeed these are the most popular Krylov methods in this area.
The runs are done on one node of a Linux cluster equipped
with a quad core Intel CPU at 2.8GHz and 16 GB of physical RAM. 
using a Portland Group Fortran~90 compiler version~9.

In Table~\ref{tab:solvers1}, we show the number of iterations required by Krylov methods to reduce 
the initial residual to $\op(10^{-5})$ starting from the zero vector. The right-hand side of the 
linear system is set up so that the initial solution is the vector of all ones.
We carry out the M-V product at each iteration using 
dense linear algebra packages, \ie~the \texttt{ZGEMV} routine available in the LAPACK library 
and we do not use preconditioning.
We observe again the remarkable effectiveness of the CORS method, that is the fastest non-Hermitian solver with 
respect to CPU time on most selected examples except GMRES with large restart.
Indeed, unrestarted GMRES may outperform all other Krylov methods and should be used when memory 
is not a concern. However reorthogonalization costs may penalize the GMRES convergence in large-scale 
applications, so using high values of restart may not be  convenient (or even not affordable for the 
memory) as shown in earlier studies~\cite{cdgs:05}. In Table~\ref{tab:solvers1} we select a value of 
50 for the restart parameter.

\begin{figure}[!ht]
 \centering
 \subfigure[On Example~2]
 { \label{fig:conv2}
   \includegraphics[width=0.52\textwidth]{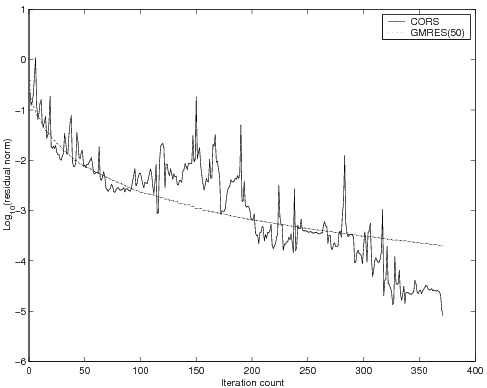}
 }
 \caption{\label{fig:convhist}Convergence histories for CORS and GMRES.}
\end{figure}

The good efficiency of CORS is even more evident on the two realistic aircraft 
problems i.e. Examples~3-4 which are very difficult to solve by iterative 
methods as no convergence is obtained without preconditioning in 3000 iterates. In Table~\ref{tab:aircrafts} 
we report on the number of iterations and on the CPU time to reduce the initial residual to $\op(10^{-3})$. 
This tolerance may be considered accurate enough for engineering purposes. In~\cite{cdgs:05} it has 
been shown that a level of accuracy of $\op(10^{-3})$ may enable a correct reconstruction of the radar 
cross section of the object. Again, CORS is more efficient than restarted GMRES on these two tough problems.
The BiCOR method also shows fast convergence. 
However, a nice feature of CORS over BiCOR is that it does not require matrix multiplications by $A^H$
on complex systems. 
This may represent an advantage when MLFMA is used because the Hermitian product often requires a 
specific algorithmic implementation~\cite{sylv:02}. In Figures~\ref{fig:convhist} we illustrate the 
convergence history of CORS and GMRES(50) on Examples~2 to show the different numerical behavior 
of the two families of solvers. The residual reduction is much smoother for GMRES along the iterations. 
We observe that in our experiments, BiCGSTAB and unsymmetric QMR algorithms 
generally converge more slowly than BiCOR and CORS. 

\begin{table}
  \centering
  \begin{small}
\begin{tabular}{clcccccc}

   Example & Description           & Size  & Memory (Gb) & Frequency (MHz) & $\kappa_1(A)$   & Geometry\\
   \hline
         1 & Sphere                & 12000 & 4.6 & 535 & $6 \cdot \op(10^5)$ & Fig.~\ref{fig:ex5} \\
         2 & Satellite             & 1699  & 0.1 & 57  & $1 \cdot \op(10^5)$ & Fig.~\ref{fig:ex6} \\
         3 & Jet prototype         & 7924  & 2.0 & 615 & $1 \cdot \op(10^7)$ & Fig.~\ref{fig:ex7} \\
         4 & Airbus A318 prototype & 23676 & 18.0 & 800 & $1 \cdot \op(10^7)$ & Fig.~\ref{fig:ex8} \\

\end{tabular}
  \caption{Characteristics of the model problems.}\label{tab:problems}
  \end{small}
\end{table}

\begin{table}
  \centering
\begin{small}
\begin{tabular}{l|cc}

  ${\raise0.7ex\hbox{Solver} \!\mathord{\left/
 {\vphantom {{Solver} {Example}}}\right.\kern-\nulldelimiterspace}
\!\lower0.7ex\hbox{Example}}$   & 1           & 2         \\
\hline

      CORS & 380 (\textbf{211}) & 371 (\textbf{11}) \\
     BiCOR & 441 (251)          & 431 (15)     \\
 GMRES(50) & $>$ 3000 ($>$ 844) & 871 (17)     \\
       QMR & 615 (508)          & 452 (24)     \\
     TFQMR & 399 (435)          & 373 (27)     \\
  BiCGSTAB & 764 (418)          & 566 (18)     \\

\end{tabular}
  \caption{Number of iterations and CPU time (in seconds) required by Krylov 
methods to reduce the initial residual to $\op(10^{-5})$.}\label{tab:solvers1}
\end{small}
\end{table}

\begin{table}
  \centering
\begin{small}
\begin{tabular}{c|cc}

${\raise0.7ex\hbox{Example} \!\mathord{\left/
 {\vphantom {{Solver} {Example}}}\right.\kern-\nulldelimiterspace}
\!\lower0.7ex\hbox{Solver}}$  &       CORS &    GMRES(50) \\
\hline

         3 & 1286 (\textbf{981}) & $>$3000 ($>$1147) \\
         4 & 924 (\textbf{5493}) & 2792 (8645) \\

\end{tabular}
\end{small}
  \caption{Number of iterations and CPU time (in seconds) for CORS and GMRES(50) to reduce the initial residual to $\op(10^{-3})$ on the two aircraft problems~\ref{fig:ex7} and~\ref{fig:ex8}. 
These problems do not converge in 3000 iterations.}\label{tab:aircrafts}
\end{table}

The large spectrum of real and complex linear systems, in both sparse and 
dense matrix computation, reported in this study illustrate the favourable 
numerical properties of the proposed unsymmetric variant of the 
Lanczos procedure for linear systems. The results indicate that our 
computational techniques are capable to solve very efficiently a large 
variety of problems. Therefore, they can be a suitable numerical tool to 
use in applications.

\section*{Acknowledgment}

We gratefully thank the EMC Team at CERFACS in Toulouse and to EADS-CCR (European Aeronautic Defence and Space - Company Corporate Research Center) in Toulouse,
for providing us with some test examples used in the numerical experiments.  This research was
supported by 973 Program (2007CB311002), NSFC (60973015) and the Project of National Defense Key Lab. (9140C6902030906).

\bibliography{mybib}

\end{document}